\documentclass[12pt]{amsart}
\usepackage{amssymb,latexsym,amsmath}
\newtheorem{prop}{Proposition}
\newcommand{\defn}{{\bf Definition: }}
\newcommand{\Notation}{{\bf Notation: }}

\newcommand{\Aside}{{\it Aside: }}
\newcommand{\EndAside}{{\it EndAside. }\smallskip}
\newcommand{\Question}{{\it Question: }}
\newcommand{\Claim}{{\it Claim: }}

\newcommand{\complex}{\mathbb{C}}
\newcommand{\matrices}{\complex^{n\times n}}
\newcommand{\sgn}{\mathrm{sgn}}
\newcommand{\Even}{\mathrm{Even}}
\newcommand{\Odd}{\mathrm{Odd}}
\newcommand{\ESym}{\mathrm{E\text{-}Sym}}
\newcommand{\ESkew}{\mathrm{E\text{-}Skew}}
\newcommand{\Circ} {\mathrm{Circ }}
\newcommand{\SCirc}{\mathrm{SCirc}}
\newcommand{\Diag}{\mathrm{Diag}}
\newcommand{\Tridiag}{\mathrm{Tridiag}}
\newcommand{\subdiag}{\mathrm{subdiag}}
\newcommand{\diag}{\mathrm{diag}}
\newcommand{\supdiag}{\mathrm{supdiag}}

\newcommand{\behnetal}{Behn, et al(2011)}
\newcommand{\bdh}{Behn, Driessel and Hentzel(2011)}
\newcommand{\catral}{Catral, et al(2009)}
\newcommand{\davis}{Davis(1979)}
\newcommand{\drew}{Drew, et al(2000)}
\newcommand{\elsner}{Elsner, et al(2003)}
\newcommand{\halmos}{Halmos(1958)}
\newcommand{\kkma}{Kailath, Kung and Morf(1979a)}
\newcommand{\kkmb}{Kailath, Kung and Morf(1979b)}
\newcommand{\ks}{Kailath and Sayed(1999)}

\newcommand{\trench}{Trench(2004)}

\begin{document}

\title
[Special Near-Toeplitz Matrices]
{A Relation Between Some Special\\
Centro-Skew, Near-Toeplitz, \\
Tridiagonal Matrices\\
And Circulant Matrices}

\date{\today} 

\keywords{Tridiagonal, Toeplitz, eigenvalue, eigenvector, centro symmetric, centro skew symmetric, sign pattern}

\subjclass{Primary: 15B05 Toeplitz, Cauchy, and related matrices; Secondary: 15B35 Sign pattern matrices; 15A18 Eigenvalues, singular values, and eigenvectors}

\maketitle

\centerline{by}

\centerline{Kenneth R. Driessel}
\centerline{Mathematics Department}
\centerline{Iowa State Universiy}

\dedicatory{}

\begin{abstract}
Let $n\ge 2$ be an integer. Let $R_n$ denote the $n\times n$ tridiagonal matrix with $-1$'s on the sub-diagonal, $1$'s on the super-diagonal, $-1$ in the $(1,1)$ entry, $1$ in the $(n,n)$ entry and zeros elsewhere. 
This paper shows that $R_n$ is closely related to a certain circulant matrix and a certain skew-circulant matrix. More precisely, let $E_n$ denote the exchange matrix which is defined by 
$E_n(i,j):=\delta(i+j,n+1)$. 
Let $E_+$ (respectively, $E_-$) be the projection defined by 
$x\mapsto (1/2)(x + E_n x)$
(respectively,
$x\mapsto (1/2)(x - E_n x)$).
Then 
$$
R_n = (\pi_n - \pi_n^T) E_+
+ (\eta_n - \eta_n^T) E_- ,
$$
where $\pi_n$ is the basic $n\times n$ circulant matrix and $\eta_n$ is the basic $n\times n$ skew-circulant matrix.
In other words, if $x$ is a vector in the range of $E_+$ then 
$R_n x = (\pi_n - \pi_n^T)x$
and if $x$ is in the range of $E_-$ then
$R_n x = (\eta_n - \eta_n^T)x$.
\end{abstract}  

\section*{Table of contents}

\begin{itemize}
\item
Introduction
\item
Circulants and Skew-Circulants
\item
Centro-Symmetric and Centro-Skew Matrices 
\item
The Relation
\item
Acknowledgements
\item
Appendix: Circulants and Skew-Circulants
\item
References
\end{itemize}


\section*{Introduction} 

We consider a special n-by-n, near-Toeplitz, tridiagonal matrix with entries from the set $\{0,1,-1\}$. In particular, we consider the tridiagonal matrix having the following form:
$$
R_n := \Tridiag (\subdiag,\diag,\supdiag)
$$
where
\begin{itemize}
\item
$\subdiag:=(-1,-1,\dots,-1)$,
\item
$\diag:=(-1,0,0,\dots,0,0,1)$, and
\item
$\supdiag:=(1,1,\dots,1,1)$.
\end{itemize}
In other words, $R_n$ is the tridiagonal matrix that 
has all $-1$'s on the subdiagonal, all 0's on the diagonal except for a $-1$ in the (1,1) entry and a 1 in 
the (n,n) entry, and has all 1's on the superdiagonal. For example, when
$n:= 5$ we have:
$$
R_5:=
\begin{pmatrix}
-1 &1  &0  &0  &0 \\
-1 &0  &1  &0  &0 \\
0  &-1 &0  &1  &0 \\
0  &0  &-1 &0  &1 \\
0  &0  &0  &-1 &1 \\
\end{pmatrix}.
$$

Here is some motivation for studying the matrices $R_n$. I say that an $n\times n$ matrix $S$ is a {\bf sign pattern matrix} if its entries come from the set $\{-1,0,+1\}$. I say that an $n\times n$ real matrix $A$ {\bf has sign pattern} $S$ if, for all $i,j=1,\dots,n$, the sign of $a_{ij}$ is the same as the sign of $s_{ij}$; in symbols,
$\sgn(a_{ij})=s_{ij}$ where $\sgn$ denotes the signum function.

\drew\ considered the sign pattern $R_n$. (They used the notation $T_n$ instead of $R_n$.) They raised the following question (among others):

\Question Does there exist a nilpotent matrix with sign pattern $R_n$ for all values of $n$?

They conjectured that the answer is yes. They showed that the answer is yes for $2\le n\le 7$. \elsner\ showed that the answer is yes for $8\le n\le 16$. (See also \catral.)

Recently, \behnetal\ proved that the answer is yes for all $n$. In particular, they proved the following result.

\begin{prop}{\bf Nilpotent sign pattern.}
For $k=1,\dots,n$, let 
$f_k:= 1/(2\sin(\theta_k))$ where 
$\theta_k:= (2k-1)\pi/(2n)$. 
Then the matrix 
$\Diag(f_1,\dots,f_n)R_n$ is nilpotent.
\end{prop}

I believe that a better understanding of the matrix $R_n$ will lead to a better understanding of this result. Recently \bdh\ solved the eigen-problem for $R_n$ in order to better understand these matrices.

The {\bf exchange matrix} $E_n$ is the $n\times n$ matrix defined by
$E_n(i,j):=\delta(i+j,n+1)$. 
An $n\times n$ matrix $H$ is {\bf centro-symmetric} if $E_n H E_n = H$. 
An $n\times n$ matrix $K$ is {\bf centro-skew} if $E_n K E_n = -K$. 
We review the theory of centro-symmetry below.
Note that the matrix $R_n$ of interest is centro-skew.

I say that a vector $x\in\complex^n$ is {\bf even} (respectively, {\bf odd}) if $E_n x = x$ (rexpectively, $E_n x = -x$). Let $\Even_n$ (respectively, $\Odd_n$) denote the subspace of $\complex^n$ consisting of the even (respectively, odd) vectors. It is easy to see that $\complex^n$  is the direct sum of the subspaces $\Even_n$ and $\Odd_n$; in symbols, $\complex^n = \Even_n \oplus \Odd_n$. 
Let $E_+$ and $E_-$ be the linear maps defined by 
\begin{align*}
E_+ := \complex^n \to \complex^n:
x\mapsto (1/2)(x+E_n x) \\
E_- := \complex^n \to \complex^n:
x\mapsto (1/2)(x-E_n x). 
\end{align*}
It is easy to see that $E_+$ (respectively, $E_-$) is the projection of $\complex^n$ onto the subspace of even (respectively, odd) vectors.

Let $\pi_n$ denote the $n\times n$ basic circulant matrix and let $\eta_n$ denote the $n\times n$ basic skew-circulant matrix. (See the section on ``Circulants and Skew-Circulants'' for the definitions.) For example, when $n=5$, 
$$
\pi_5:=
\begin{pmatrix}
0 &1 &0 &0 &0 \\
0 &0 &1 &0 &0 \\
0 &0 &0 &1 &0 \\
0 &0 &0 &0 &1 \\
1 &0 &0 &0 &0 \\
\end{pmatrix}
$$
and
$$
\eta_5:=
\begin{pmatrix}
0  &1 &0 &0 &0 \\
0  &0 &1 &0 &0 \\
0  &0 &0 &1 &0 \\
0  &0 &0 &0 &1 \\
-1 &0 &0 &0 &0 \\
\end{pmatrix}.
$$

In this paper I prove the following result. 
It says that the restriction of $R_n$ to the subspace of even vectors is a circulant matrix and its restriction to the subspace of odd vectors is a skew-circulant matrix.

\begin{prop} {\bf The Relation.}
The special tridiagonal matrix $R_n$, the basic circulant $\pi_n$ and the basic skew-circulant $\eta_n$ satisfy
$$
R_n =(\pi_n - \pi_n^T) E_+ + (\eta_n - \eta_n^T)E_- .
$$
In other words, for all $x\in \complex^n$,
\begin{itemize}
\item
if $x$ is even then $R_n x = (\pi_n - \pi_n^T)x$, and
\item
if $x$ is odd then $R_n x = (\eta_n - \eta_n^T)x$.
\end{itemize}
\end{prop}

The following calculations illustrate this result when $n=5$:
$$
\begin{pmatrix}
-1 &1  &0  &0  &0 \\
-1 &0  &1  &0  &0 \\
0  &-1 &0  &1  &0 \\
0  &0  &-1 &0  &1 \\
0  &0  &0  &-1 &1 \\
\end{pmatrix} 
\begin{pmatrix}
x_1  \\
x_2  \\
x_3  \\
x_2  \\
x_1 \\
\end{pmatrix}
=
\begin{pmatrix}
x_2 - x_1 \\
x_3 - x_1 \\
0  \\
x_1 - x_3 \\
x_1 - x_2 \\
\end{pmatrix},
$$
$$
\begin{pmatrix}
0  &1 &0 &0 &-1 \\
-1  &0 &1 &0 &0 \\
0  &-1 &0 &1 &0 \\
0  &0 &-1 &0 &1 \\
1 &0 &0 &-1 &0 \\
\end{pmatrix}
\begin{pmatrix}
x_1  \\
x_2  \\
x_3  \\
x_2  \\
x_1 \\
\end{pmatrix}
=
\begin{pmatrix}
x_2 - x_1 \\
x_3 - x_1 \\
0  \\
x_1 - x_3 \\
x_1 - x_2 \\
\end{pmatrix},
$$
$$
\begin{pmatrix}
-1 &1  &0  &0  &0 \\
-1 &0  &1  &0  &0 \\
0  &-1 &0  &1  &0 \\
0  &0  &-1 &0  &1 \\
0  &0  &0  &-1 &1 \\
\end{pmatrix}
\begin{pmatrix}
x_1  \\
x_2  \\
0  \\
-x_2  \\
-x_1 \\
\end{pmatrix}
=
\begin{pmatrix}
x_2 - x_1 \\
- x_1 \\
-2 x_2  \\
- x_1  \\
x_2 - x_1 \\
\end{pmatrix},
$$
$$
\begin{pmatrix}
0  &1  &0  &0  &1 \\
-1 &0  &1  &0  &0 \\
0  &-1 &0  &1  &0 \\
0  &0  &-1 &0  &1 \\
-1 &0  &0  &-1 &0 \\
\end{pmatrix}
\begin{pmatrix}
x_1  \\
x_2  \\
0  \\
-x_2  \\
-x_1 \\
\end{pmatrix}
=
\begin{pmatrix}
x_2 - x_1 \\
- x_1 \\
-2 x_2  \\
- x_1  \\
x_2 - x_1 \\
\end{pmatrix}.
$$

Here is a summary of the contents. In the section with title ``Circulants and Skew-Circulants'', I present a brief review of the theory of these classes of matrices. (In the appendix with the same title, I present a more extensive review.) In the section with title ``Centro-Symmetric and Centro-Skew Matrices'', I present a review of the theory of these classes of matrices. These two sections are included for the reader's convenience. In the section with title ``The 
Relation'', I prove that the relation between 
$R_n$, $\pi_n$, and $\eta_n$
described above holds.

In this paper I repeat some definitions to accommodate ``grasshopper'' readers.


\section*{Circulants and Skew-Circulants} 

In this section we review the theory of circulant and skew-circulant matrices (which are defined below). For more details (including proofs), see the appendix on circulants and skew-circulants. 

\defn Let $C$ be an $n\times n$ matrix. Then $C$ is a {\bf circulant matrix} or {\bf circulant} for short, if there exist scalars $c_1,\dots,c_n$, such that $C = \Circ(c_1,c_2,\dots,c_n)$ where
\begin{equation*}
\Circ(c_1,c_2,\dots,c_n):=
\begin{pmatrix}
c_1     & c_2   & c_3    & \dots & c_n \\
c_n     &c_1   &c_2    &\dots &c_{n-1} \\
c_{n-1} &c_n   &c_1    &\dots &c_{n-2} \\
\vdots  &\vdots &\vdots &\ddots &\vdots \\
c_2 &c_3 &c_4 &\dots &c_1 
\end{pmatrix}.
\end{equation*}
The $n\times n$ matrix $\pi:=\pi_n:=\Circ(0,1,0,\dots,0)$ is the {\bf basic circulant} with size $n$. (I also use the symbol $\pi$ to with its usual meaning - namely, the ratio of the circumference of a circle to its diameter. I assume that the reader can distinguish these two meanings from the context.)

\begin{prop}
The basic $n\times n$ circulant satisfies the following equations:
$
\pi^n = I,\quad \pi^T = \pi^\star = \pi^{-1} = \pi^{n-1}.
$
The minimum polynomial of $\pi$ is $\lambda^n-1$.
For every sequence of scalars $c_1,c_2,\dots,c_n$,
$$
\Circ(c_1,c_2,\dots,c_n)
= c_1 I + c_2 \pi + c_3 \pi^2 + \cdots c_n \pi^{n-1}.
$$

\end{prop}

\Notation For any $a:= (a_1,a_2,\dots, a_n)$, let 
$p_a$ denote the polynomial defined by
$p_a(t):= a_1 + a_2 t +\cdots +a_n t^{n-1}$. 

From the last proposition, we see that we can write 
$\Circ(c) = p_c(\pi)$ where $c:=(c_1,\dots,c_n)$. 

\Notation Let $n$ be a positive integer and let $
\omega:=\omega_n := \exp(i\ 2\pi/n)$. 

\defn The $n\times n$ matrix $F$, defined by $F^\star (i,j):= \omega^{(i-1)(j-1)}$
is the {\bf Fourier matrix of order $n$}.

\begin{prop}
The Fourier matrix satisfies the following equations:
$$
F = F^T, \quad F^\star = (F^\star)^T = \bar{F},\quad
F = \bar{F}^\star. 
$$
\end{prop}

\begin{prop}
The Fourier matrix is unitary.
\end{prop}

\Notation Let $\Omega:=\Omega_n:=\Diag(1,\omega,\omega^2,\dots,\omega^{n-1})$
denote the diagonal matrix with diagonal entries $1,\omega,\dots,\omega^{n-1}$. 

\begin{prop}{\bf Spectral decomposition of the basic circulant.}
The basic circulant $\pi$ satifies the equation $\pi = F^\star \Omega F$. Hence the eigen-pairs  of $\pi$ are $\omega^{k-1}$, $f_k$, 
for $k=1,\dots,n$,
where $f_k$ is the $k$th column of the matrix $F^\star$.
\end{prop}

\begin{prop}
{\bf Spectral decomposition for circulants.} Let
$c:= (c_1,\dots,c_n)$.
Then $\Circ(c) = F^\star p_c(\Omega) F$. 
Hence the eigen-pairs of $\Circ(c)$ are 
$p_c(\omega^{k-1})$, $f_k$,
for $k=1,2,\dots,n$, where $f_k$ is the $k$th column of $F^\star$. 
\end{prop}

\defn Let $S$ be an $n\times n$ matrix. Then $S$ is a {\bf skew-circulant matrix} or {\bf skew-circulant} for short, if there exist scalars $a_1,\dots,a_n$, such that $S = \SCirc(a_1,a_2,\dots,a_n)$ where
\begin{equation*}
\SCirc(a_1,a_2,\dots,a_n):=
\begin{pmatrix}
a_1     & a_2   & a_3    & \dots & a_n \\
-a_n     &a_1   &a_2    &\dots &a_{n-1} \\
-a_{n-1} &-a_n   &a_1    &\dots &a_{n-2} \\
\vdots  &\vdots &\vdots &\ddots &\vdots \\
-a_2 &-a_3 &-a_4 &\dots &a_1 
\end{pmatrix}.
\end{equation*}
The $n\times n$ matrix $\eta:=\eta_n:=\SCirc(0,1,0,\dots,0)$ is the {\bf basic skew-circulant} with size $n$. 

\begin{prop}
The basic $n\times n$ skew-circulant $\eta$ satisfies the following equations:
$
\eta^n = -I,\quad \eta^T = \eta^\star = \eta^{-1} = -\eta^{n-1}.
$
The minimum polynomial of $\eta$ is $\lambda^n+1$.
For every sequence of scalars $a_1,a_2,\dots,a_n$,
$$
\SCirc(a_1,a_2,\dots,a_n)
= a_1 I + a_2 \eta + a_3 \eta^2 + \cdots a_n \eta^{n-1}.
$$
\end{prop}

Note that, if $a:=(a_1,\dots,a_n)$ then
$\SCirc(a) = p_a(\eta)$. 

\Notation Let $n$ be a positive integer and let $\sigma:=\sigma_n := \exp(i \pi/n)$. Let
$\Omega^{1/2}:=\Omega_n^{1/2}:= \Diag(1,\sigma_n,\sigma_n^2,\dots,\sigma_n^{n-1})$.

Note that $\omega_n=\sigma_n^2$ and hence
$\Omega=(\Omega^{1/2})^2$. 

\begin{prop}{\bf The relation between the basic circulant and the basic skew-circulant.}
The basic matrices $\pi$ and $\eta$ satisfy the following equation:
$$ 
\eta=\sigma \Omega^{1/2}\pi \bar{\Omega}^{1/2}
 =\Omega^{1/2}(\sigma\pi) \bar{\Omega}^{1/2}.
$$
\end{prop}

\Notation Let the $n\times n$ matrix $H:=H_n$, be defined by 
$$H^\star := \Omega^{1/2} F^\star.$$

\begin{prop}
The matrix $H$ is unitary.
\end{prop}

\begin{prop}{\bf Spectral decomposition of the basic skew-circulant.}
The basic skew-circulant $\eta$ satifies the eauation $\eta = H^\star (\sigma\Omega) H$. 
Hence the eigen-pairs of $\eta$ are $\sigma^{2k-1}$, 
$h_k$, for $k=1,\dots,n$, where $h_k$ is the $k$th column of the matrix $H^\star$. 
\end{prop}

\begin{prop}
{\bf Spetral decomposition for skew-circulants.} 
Let $a:=(a_1,\dots,a_n)$. Then
$\SCirc(a)=H^\star p_a(\sigma\Omega)H$. 
Hence the eigen-pairs of $\Circ(a)$ are 
$p_a(\sigma^{2k-1})$, $h_k$,
for $k=1,2,\dots,n$, where $h_k$ is the $k$th column of the matrix $H^\star$.
\end{prop}


\section*{Centro-Symmetric and Centro-Skew Matrices} 

In this section, we review the (fairly well-known) theory of centro-symmetric and centro-skew-symmetric matrices. (These terms are defined below).  For more details and generalizations see \trench\ and the references in that paper. 

\defn The {\bf exchange} (or {\bf flip}) {\bf matrix} $E_n$ is defined by  
$$
E_n(i,j) := \delta(i+j,n+1)
$$
where $\delta$ is the Kronecker delta. Here is a picture of $E_4$:
$$
E_4:=
\begin{pmatrix}
0 &0 &0 &1 \\
0 &0 &1 &0 \\
0 &1 &0 &0 \\
1 &0 &0 &0
\end{pmatrix}
$$
The exchange matrix has ones on the ``per-diagonal'' (or ``counter-diagonal'') and zeros elsewhere. I often drop the subscript; in other words, I often write $E$ instead of $E_n$. 

Note that every exchange matrix $E$ is a symmetric permutation matrix which satisfies $E^2=I$. (Recall that a linear transformation $L$ on a vector space is an {\bf involution} if $L^2=I$.) It follows that the eigenvalues of $E$ are $1$ and $-1$. 

\defn Let $x$ be a vector in $\complex^n$. If $Ex=x$ then $x$ is an {\bf even vector}; if $Ex=-x$ then $x$ is an {\bf odd vector}. 

For example, if $n=5$, even vectors have the form $(x_1,x_2,x_3,x_2,x_1)^T$ and odd vectors have the form  $(x_1,x_2,0,-x_2,-x_1)^T$. I shall use $\Even_n$ (or simply $\Even$) to denote the space of even vectors and $\Odd_n$ (or simply $\Odd$) to denote the space of odd vectors. Note that $\complex^n$ is the direct sum of the even and odd subspaces; in symbols, $\complex^n=\Even\oplus\Odd$.

\Notation Let the linear maps $E_+$ and $E_-$ be defined as follows:
\begin{align*}
E_+ &:= \complex^n \to \complex^n: x \mapsto (1/2)(x+Ex) \\
E_- &:= \complex^n \to \complex^n: x \mapsto (1/2)(x-Ex).
\end{align*}

The following result lists some of the elementary properties of these linear maps. The proofs are straight-forward. 

\begin{prop}
The linear maps $E_+$ and $E_-$ satisfy the following conditions:
\begin{itemize}
\item
The range of $E_+$ equals the subspace of even vectors; the range of $E_-$ equals the subspace of odd vectors.
\item
$E = E_+ - E_-,\ I = E_+ + E_-.$
\item
$ E_+^2 = E_+,\ E_-^2 = E_-,\  
E_+ E_- = E_- E_+ = 0.$
\item
$E E_+ = E_+ E = E_+,\ 
E E_- = E_- E = -E_-.$
\end{itemize}
\end{prop}

\begin{proof}
The equations in the last item follow from the previous ones: in particular, we have
\begin{align*}
EE_+ &= (E_+ - E_-)E_+ = E_+^2 = E_+,\\
EE_- &= (E_+ - E_-)E_- = -E_-^2 = -E_-.
\end{align*}
\end{proof}

\Aside Note that $E=E_+ - E_-$ is the ``spectral form'' of $E$. See \halmos.
\EndAside

\Notation Define the linear map $\phi:=\phi_n$ as follows:
$$
\phi_n := \matrices\to\matrices:X\mapsto EXE.
$$
Note that this map is an involution. 

\defn A linear map $A:\complex^n\to\complex^n$ is  
{\bf centro-symmetric} (or
{\bf $E$-symmetric}) if $\phi_n(A)=A$; 
a linear map $K:\complex^n\to\complex^n$ is  
{\bf centro-skew-symmetric} (or 
{\bf centro-skew} for short, or  
{\bf $E$-skew} ) if $\phi_n(K)=-K$. 

I shall use $\ESym(n)$ (or simply $\ESym$) to denote the subspace of centro-symmetric matrices and I shall use $\ESkew(n)$ (or simply $\ESkew$) to denote the subspace of centro-skew matrices.  

When $n=4$, centro-symmetric matrices have the form
$$
\begin{pmatrix}
a_1 &a_2 &a_3 &a_4 \\
a_5 &a_6 &a_7 &a_8 \\
a_8 &a_7 &a_6 &a_5 \\
a_4 &a_3 &a_2 &a_1
\end{pmatrix}
$$
and centro-skew matrices have the form
$$
\begin{pmatrix}
k_1 &k_2 &k_3 &k_4 \\
k_5 &k_6 &k_7 &k_8 \\
-k_8 &-k_7 &-k_6 &-k_5 \\
-k_4 &-k_3 &-k_2 &-k_1
\end{pmatrix}.
$$

The following result lists some of the basic properties of centro-symmetric and centro-skew matrices. I omit the striaght-forward proofs.

\begin{prop}{\bf Basic properties of centro-symmetry.}
\begin{itemize}
\item
The space $\matrices$ of matrices is the direct sum of the subspace of centro-symmetric matrices and the subspace of centro-skew marices; in symbols,
$\matrices=\ESym\oplus\ESkew$.
\item 
If $A_1$ and $A_2$ are centro-symmetric then the $A_1 A_2$ is centro-symmetric; if $A$ is centro-symmetric and $K$ is centro-skew then $AK$ and $KA$ are centro-skew; if $K_1$ and $K_2$ are centro-skew then $K_1K_2$ is centro-symmetric.
\item
Let $A$ be centro-symmetric. If the vector $x$ is even then the vector $Ax$ is even; if the vector $x$ is odd then the vector $Ax$ is odd.
\item
Let $K$ be centro-skew. If the vector $x$ is even then the vector $Kx$ is odd; if the vector $x$ is odd then the vector $Kx$ is even. 
\end{itemize}
\end{prop}

Here is the ``multiplication table'' for centro-symmetric and centro-skew matrices:
\smallskip

\begin{tabular}{c | c | c}
$\ast$   & $\ESym$  & $\ESkew$ \\ \hline
$\ESym$  & $\ESym$  & $\ESkew$ \\ \hline
$\ESkew$ & $\ESkew$ & $\ESym$ 
\end{tabular}
\smallskip

\begin{prop}{\bf Decomposition of linear maps.}
\begin{itemize}\item
Every linear map $L:\complex^n\to\complex^n$ can be decomposed as follows:
$$
L = E_+ L E_+ + E_+ L E_- + E_- L E_+ + E_- L E_-.
$$
\item
Let $A:\complex^n\to\complex^n$ be a linear map. Then $A$ is $E$-symmetric iff 
$$
A = E_+ A E_+ + E_- A E_-.
$$
\item
Let $K:\complex^n\to\complex^n$ be a linear map. Then $K$ is $E$-skew iff 
$$
K = E_+ K E_- + E_- K E_+.
$$
\end{itemize}
\end{prop}

\begin{proof}
The first item follows from the equation
$$
L=(E_+ + E_-)L(E_+ + E_-).
$$

Now assume that $EAE=A$. Then
\begin{align*}
A &= EAE = (E_+ - E_-)A(E_+ - E_-) \\
&= E_+ A E_+ - E_+AE_-
- E_- A E_+ + E_-AE_-.
\end{align*}
Also, by the first item, we have
$$
A = E_+ A E_+ + E_+AE_-
  + E_- A E_+ + E_-AE_-.
$$
It follows that 
$E_+AE_- + E_-AE_+ = 0.$

Now assume that 
$A = E_+AE_+ + E_-AE_- .$
Then
\begin{align*}
EAE &= E(E_+AE_+ + E_-AE_-)E 
= EE_+AE_+E + EE_-AE_-E \\
&= E_+AE_+ + E_-AE_- = A. 
\end{align*}

Now assume that $EKE=-K$. Then
\begin{align*}
-K &= EKE = (E_+ - E_-)K(E_+ - E_-) \\
&= E_+KE_+ - E_+KE_-
- E_-KE_+ + E_-KE_-.
\end{align*}
Also, by the first item, we have
$$
-K = -(E_+KE_+ + E_+KE_-
     + E_-KE_+ + E_-KE_-).
$$
It follows that 
$E_+KE_+ + E_-KE_- = 0.$

Now assume that 
$K = E_+KE_- + E_-KE_+ .$
Then
\begin{align*}
EKE &= E(E_+KE_- + E_-KE_+)E 
= EE_+KE_-E + EE_-KE_+E \\
&= -E_+KE_- - E_-KE_+ = -K. 
\end{align*}
\end{proof}

\begin{prop}{\bf Decomposition of solutions.}
\begin{itemize}
\item
Let $A$ be $E$-symmetric. Then $Az=w$ iff
$AE_+z = E_+ w$ and $AE_-z = E_-w$.
\item
Let $K$ be $E$-skew. Then $Kz=w$ iff
$KE_+z = E_- w$ and $AE_-z = E_+w$.
\end{itemize}
\end{prop}

\begin{proof}
Assume $AE_+z = E_+w$ and $AE_- z = E_-w$.
Then 
$$
Az = A(E_+ + E_-)z 
= AE_+ z + AE_- z
= E_+ w + E_- w = w.
$$

Assume $Az=w$. Then 
$AE_+z + AE_- z = E_+ w + E_- w.$
Note
$E(AE_+ z) = AEE_+ z = AE_+ z,$ and
$E(AE_- z) = AEE_- z = -AE_+ z.$
Hence $AE_+ z$ is even and $AE_- z$ is odd. 
Since $\complex^n$ is the direct sum of the even subspace and the odd subspace, we get
$AE_+ z = E_+ w$ and $AE_- z = E_- w$.

Assume $KE_+z = E_-w$ and $KE_- z = E_+w$.
Then 
$$
Kz = K(E_+ + E_-)z 
= KE_+ z + KE_- z
= E_- w + E_+ w = w.
$$

Assume $Kz=w$. Then 
$KE_+z + KE_- z = E_+ w + E_- w.$
Note
$E(KE_+ z) = -KEE_+ z = -KE_+ z,$ and 
$E(KE_- z) = -KEE_- z =  KE_- z.$
Hence $KE_+ z$ is odd and $KE_- z$ is even. 
Since $\complex^n$ is the direct sum of the even subspace and the odd subspace, we get
$KE_+ z = E_- w$ and $KE_- z = E_+ w$
\end{proof}

\begin{prop}{\bf Decomposition of eigenvectors for $E$-symmetric maps.}
Let $A$ be $E$-symmetric. Then
$Az=\lambda z$ iff 
$AE_+ z = \lambda E_+ z$ and $AE_- z = \lambda E_- z$.
Furthermore, if $Az=\lambda z$ and $E_+ z\neq 0$ then $E_+z$ is an even eigenvector of $E_+AE_+$ associated with $\lambda$ and, if $E_-z\neq0$ then 
$E_-z$ is an odd eigenvector of $E_-AE_-$ associated with $\lambda$.  
\end{prop}

\begin{proof}
From the result on decomposition of solutions we have 
$Az=\lambda z$ iff
$AE_+z = \lambda E_+z$ and 
$AE_-z = \lambda E_-z$.

The proof of the last assertion in the proposition is straight-forward.
\end{proof}

\begin{prop}{\bf Decomposition of eigenvectors for $E$-skew maps.}
Let $K$ be $E$-skew.
\begin{itemize}
\item
Then
$Kz=\lambda z$ iff
$KE_- z = \lambda E_+ z$ and 
$KE_+ z = \lambda E_- z$.
\item 
If $(\lambda,z)$ is an eigen-pair for $K$ then so is
$(-\lambda,Ez)$.
\item
If $(\lambda,z)$ is an eigen-pair for $K$ and $\lambda\neq 0$ then $z$ is neither even nor odd; in other words, $z\notin \Even\cup \Odd$. 
\item
If $K$ is singular, then the null space of $K$ has a basis in $\Even\cup\Odd$.
\end{itemize}
\end{prop}

\begin{proof}
From the result on decomposition of solutions we have $Kz=\lambda z$ iff $KE_+z=\lambda E_-z$ and
$KE_-z=\lambda E_+z$.

Now assume $Kz=\lambda z$. Then
$$
K(Ez)=-(EKE)Ez = -EKz = -\lambda Ez.
$$

Now assume $\lambda\neq 0$. Then $z$ and $Ez$ are independent because they are associated with different eigenvalues. In particular, $Ez\neq z$ and $Ez\neq -z$. Hence $z$ is neither even nor odd. 

Now assume that $K$ is singular. By the second item, we have $Kz=0$ implies $KEz=0$. Hence
$KE_+z=0$ and $KE_-z=0$. 
Let $\{z_1,\dots,z_k\}$ be a basis for the null space of $K$. Then the set
$$
\{E_+z_1,\dots,E_+z_k,E_-z_1,\dots,E_-z_k\}
$$
spans the null space and hence we can chose a basis from this set. 
\end{proof}

\Aside
The decomposition results given above can also be described in terms of block matrices. In particular,\trench\ describes the results that way.  

Let $p_1,\dots,p_r$ be an orthonormal basis for $\Even_n$ and let $q_1,\dots,q_s$ be a basis for $\Odd_n$. Note that $n=r+s$. Let 
$P:=(p_1,\dots,p_r)$ be the $n\times r$ matrix with the $p_j$ as columns and let 
$Q:=(q_1,\dots,q_s)$ be the $n\times s$ matrix with the $q_j$ as columns. 
The following result lists some of the elementary properties of the matrices $P$ and $Q$. I omit the easy proof. 

\begin{prop}
The matrices $P$ and $Q$ defined above satisfy the following conditions.
\begin{itemize}
\item
The matrix $PP^\star$ is a projection of $\complex^n$ onto $\Even$; the matrix $QQ^\star$ is a projection of $\complex^n$ onto $\Odd$.
\item
$ E = PP^\star - QQ^\star$,
$I = PP^\star + QQ^\star$.
\item
The matrix $P^\star P$ equals the identity matrix on the even subspace and the matrix $Q^\star Q$ equals the identity matrix on the odd subspace. 
\item
$P^\star Q = 0$ and $Q^\star P = 0$.
\item
The matrix
$(P,Q):= (p_1,\dots,p_r,q_1,\dots,q_s)$ is unitary.
\end{itemize}
\end{prop}

We can write any matrix $X$ in $\matrices$ in block form so that it conforms to the even-odd decomposition of the space $C^n$ as follows:
$$
X = (P,Q)
\begin{pmatrix}
X_{11} & X_{12} \\
X_{21} & X_{22}
\end{pmatrix}
\begin{pmatrix}
P^\star \\
Q^\star
\end{pmatrix}
$$
where 
$X_{11}:= P^\star X P,
X_{12}:= P^\star X Q,
X_{21}:= Q^\star X P,$
and
$X_{22}:= Q^\star X Q$.

The following proposition provides a characterization of centro-symmetric and centro-skew matrices in terms of this block decomposition.  

\begin{prop}
{\bf Block characterization of centro-symmetric and centro-skew.} Let $P$ and $Q$ be defined as above.
\begin{itemize}
\item
For every matrix $A$ in $\matrices$, $A$ is centro-symmetric iff
$$
A = (P,Q)
\begin{pmatrix}
P^\star A P & 0 \\
0        & Q^\star A Q
\end{pmatrix}
\begin{pmatrix}
P^\star \\
Q^\star
\end{pmatrix}.
$$
\item
For every matrix $K$ in $\matrices$, $K$ is centro-skew iff
$$
K = (P,Q)
\begin{pmatrix}
0 & P^\star K Q \\
Q^\star K P & 0
\end{pmatrix}
\begin{pmatrix}
P^\star \\
Q^\star
\end{pmatrix}.
$$
\end{itemize}
\end{prop}

The conclusions of the two parts of the proposition can also be written as follows:
\begin{align*}
A &= PP^\star A PP^\star + QQ^\star A QQ^\star, \\
K &= QQ^\star K PP^\star + PP^\star K QQ^\star.
\end{align*}

\EndAside


\section*{The Relation}

In this section, we prove that there is a simple relation between the special tridiagonal matrix $R_n$,  the basic circulant matrix $\pi_n$ and the basic skew-circulant matrix $\eta_n$.

Note that $R_n$ is centro-skew. 
Using the $E_+$ and $E_-$ notation from above, we have (from the result on decomposition of linear maps) that
$$
R_n = E_+ R_n E_- + E_- R_n E_+.
$$ 
We also know that centro-skew matrices map even vectors to odd ones and odd vectors to even ones. Consequently, we can write the last equation more concisely as follows:
$$
R_n = R_n E_- + R_n E_+. 
$$

\begin{prop}{\bf The relation.} The following equations hold between the special tridiagonal matrix $R_n$, the basic circulant matrix $\pi_n$ and the basic skew-circulant $\eta_n$:
$$
R_n E_- = (\eta_n - \eta_n^T)E_-\ 
\mathrm{ and }\ 
R_n E_+ = (\pi_n - \pi_n^T)E_+\  
.
$$
Hence
$$
R_n = (\eta_n - \eta_n^T) E_- + (\pi_n - \pi_n^T) E_+
$$
\end{prop}

In other words, for all $x$ in $\complex^n$,
\begin{itemize}
\item
if $x$ is even then 
$R_n x = (\pi_n - \pi_n^T)x$, and
\item
if $x$ is odd then
$R_n x = (\eta_n - \eta_n^T)x.$
\end{itemize}

\begin{proof}
We shall use descriptions of the matrices of interest that highlight the similarities and differences between them. 

Let $Z_n$ denote the {\bf lower shift matrix} which is the matrix that has 1's on the subdiagonal and 0's elsewhere. For example, when $n=4$, we have
$$
Z_4 := 
\begin{pmatrix}
0 & 0 & 0 & 0  \\
1 & 0 & 0 & 0  \\
0 & 1 & 0 & 0  \\
0 & 0 & 1 & 0 
\end{pmatrix}.
$$

\Aside
Here we follow notation and terminology used by T. Kailath. See, for example, \kkma, \kkmb\  or \ks.
\EndAside

Note that 
$$
R_n = Z_n^T - Z_n - e_1 e_1^T + e_n e_n^T
$$
where $e_k$ denotes the column vector which has 1 in the $k$th coordinate and 0's elsewhere. 
Note that 
$$
\pi_n - \pi_n^T 
= Z_n^T - Z_n - e_1 e_n^T + e_n e_1^T
$$
and
$$
\eta_n - \eta_n^T = 
Z_n^T - Z_n + e_1 e_n^T - e_n e_1^T.
$$

\Claim If $x$ is even then 
$(R_n - (\pi_n - \pi_n^T))x=0$. 

We have
\begin{align*}
R_n  - (\pi_n - \pi_n^T) 
&= (-e_1 e_1^T + e_n e_n^T) 
- (-e_1 e_n^T + e_n e_1^T) \\
&= (e_n+e_1)(e_n-e_1)^T.
\end{align*}

Note that vectors of the form $e_k+e_{n+1-k}$ span the subspace of even vectors. We now simply calculate as follows:
$$
(e_n+e_1)(e_n-e_1)^T(e_k+e_{n+1-k})
= (e_n+e_1)
(\delta_{nk} +\delta_{n(n+1-k)} 
- \delta_{1k}- \delta_{1(n+1-k)}).
$$
Now 
\[
(\delta_{nk} +\delta_{n(n+1-k})
- \delta_{1k}- \delta_{1(n+1-k)})
=
\begin{cases}
0+0-0-0, &\text{if $k\neq 1$ and $k\neq n$}\\ 
0+1-1+0, &\text{if $k=1$}\\
1+0-0-1, &\text{if $k=n$.}
\end{cases}
\]
This completes the proof of the claim.

\Claim If $x$ is odd then
$R_n x - (\eta_n - \eta_n^T)x=0.$

We have
\begin{align*}
R_n  - (\eta_n + \eta_n^T) 
&= (-e_1 e_1^T + e_n e_n^T) 
- (e_1 e_n^T - e_n e_1^T) \\
&= (e_n-e_1)(e_n+e_1)^T.
\end{align*}

Note that vectors of the form $e_k-e_{n+1-k}$ span the subspace of odd vectors. We now simply calculate as follows:
$$
(e_n+e_1)(e_n-e_1)^T(e_k-e_{n+1-k})
= (e_n+e_1)
(\delta_{nk} +\delta_{n(n+1-k)} 
- \delta_{1k}- \delta_{1(n+1-k)}).
$$
Now 
\[
(\delta_{nk} +\delta_{n(n+1-k)} 
- \delta_{1k}- \delta_{1(n+1-k)})
=
\begin{cases}
0+0-0-0, &\text{if $k\neq 1$ and $k\neq n$}\\ 
0+1-1+0, &\text{if $k=1$}\\
1+0-0-1, &\text{if $k=n$.}
\end{cases}
\]
\end{proof}


\section*{Acknowledgements} 

I discussed the ideas in this report a number of times with James Wilson (Mathematics, Iowa State University); I thank him for his comments. Irv Hentzel (Mathematics, Iowa State University) encouraged me to study sign pattern matrices; I thank him for his continued encouragement and for his constructive comments concerning this work. During his visit to Iowa, Antonio Behn (Mathematics, University of Chile) and I often discussed sign pattern matrices; I thank him for his comments on my thoughts. Wayne Barrett (Mathematics, Brigham Young University) read some of my related work ; I thank him for his comments and encouragement. During the last three decades, Tom Kailath (Engineering, Stanford University, emeritus) repeatedly encouraged me to study Toeplitz matrices; I thank him for his guidance. Wolfgang Kliemann (Chair) arranged my affiliation with the Mathematics Department at Iowa State University; I thank him.


\section*{Appendix: Circulants and Skew-Circulants} 

In this appendix I follow \davis.

\defn Let $C$ be an $n\times n$ matrix. Then $C$ is a {\bf circulant matrix} or {\bf circulant} for short, if there exist scalars $c_1,\dots,c_n$, such that $C = \Circ(c_1,c_2,\dots,c_n)$ where
\begin{equation*}
\Circ(c_1,c_2,\dots,c_n):=
\begin{pmatrix}
c_1     & c_2   & c_3    & \dots & c_n \\
c_n     &c_1   &c_2    &\dots &c_{n-1} \\
c_{n-1} &c_n   &c_1    &\dots &c_{n-2} \\
\vdots  &\vdots &\vdots &\ddots &\vdots \\
c_2 &c_3 &c_4 &\dots &c_1 
\end{pmatrix}.
\end{equation*}

Note that the $n\times n$ circulants form a linear subspace of the space of all $n\times n$ matrices.

\defn We use the notation $\pi:=\pi_n:=\Circ(0,1,0,\dots,0)$ for the {\bf basic circulant} with size $n$. 

\davis\ says the basic circulant ``plays a fundamental role in the theory of circulants''. Here is a picture of $\pi$ when $n=4$:
$$
\pi_4=
\begin{pmatrix}
0     &1   &0    &0 \\
0     &0   &1    &0 \\
0     &0   &0    &1 \\
1     &0   &0    &0 
\end{pmatrix}.
$$
Note that $\pi$ is a permutation matrix. 
The following result is easy to verify. We omit the proof.

\begin{prop}
The basic $n\times n$ circulant satisfies the following equations:
$$
\pi^n = I,\quad \pi^T = \pi^\star = \pi^{-1} = \pi^{n-1}.
$$
The minimum polynomial of $\pi$ is $\lambda^n-1$.
For every sequence of scalars $c_1,c_2,\dots,c_n$,
$$
\Circ(c_1,c_2,\dots,c_n)
= c_1 I + c_2 \pi + c_3 \pi^2 + \cdots c_n \pi^{n-1}.
$$

\end{prop}

\Notation For any $a:= (a_1,a_2,\dots, a_n)$, let 
$p_a$ denote the polynomial defined by
$p_a(t):= a_1 + a_2 t +\cdots +a_n t^{n-1}$. 

From the last proposition, we see that we can write 
$\Circ(c) = p_c(\pi)$ where $c:=(c_1,\dots,c_n)$. From this it is clear that all circulants of the same order commute. In other words, the circulants of a given order form a commutative algebra generated by the single matrix $\pi$. 

\Notation Let $n$ be a positive integer and let $
\omega:=\omega_n := \exp(i\ 2\pi/n) = \cos(2\pi/n) + i\sin(2\pi/n)$. 

\defn The $n\times n$ matrix $F$, defined by $$F^\star (i,j):= \omega^{(i-1)(j-1)}$$
is the {\bf Fourier matrix of order $n$}.

Here is a picture of $F^\star$:
$$
F^\star := \frac{1}{\sqrt n}
\begin{pmatrix}
1      &1        &1           & \dots 
&1 \\
1      &\omega   &\omega^2    &\dots 
&\omega^{n-1} \\
1      &\omega^2 &\omega^4  &\dots 
&\omega^{2(n-1)} \\
\vdots &\vdots   &\vdots      &\ddots 
&\vdots \\
1      &\omega^{n-1}&\omega^{2(n-1)} &\dots 
&\omega^{(n-1)(n-1)} 
\end{pmatrix}.
$$

The following properties of the matrix $F$ are obvious.  

\begin{prop}
The Fourier matrix satisfies the following equations:
$$
F = F^T, \quad F^\star = (F^\star)^T = \bar{F},\quad
F = \bar{F}^\star. 
$$
\end{prop}

Recall that an $n\times n$ matrix $U$ is {\bf unitary} if $UU^\star = U^\star U = I$.  

\begin{prop}
The Fourier matrix is unitary.
\end{prop}

\begin{proof}(sketch)
Recall the finite ``geometric sum identity'' which holds for any scalar $t$ different than 1:
$$
1 + t + t^2 + \cdots t^{n-1} = \frac{1-t^n}{1-t}.
$$ 
From this identity, we get
$$
1 + \omega^{2(j-k)} + \omega^{3(j-k)} + \cdots +
\omega^{(n-1)(j-k)} 
= 
\begin{cases}
n, &\text{ if $j=k$;}\\
0  &\text{ otherwise}.
\end{cases}
$$
\end{proof}

\Notation Let $\Omega:=\Omega_n:=\Diag(1,\omega,\omega^2,\dots,\omega^{n-1})$
denote the diagonal matrix with diagonal entries $1,\omega,\dots,\omega^{n-1}$. 

\begin{prop}{\bf Spectral decomposition of the basic circulant.}
The basic circulant $\pi$ satifies the equation $\pi = F^\star \Omega F$. Hence the eigen-pairs  of $\pi$ are $\omega^{k-1}$, $f_k$, 
for $k=1,\dots,n$,
where $f_k$ is the $k$th column of the matrix $F^\star$.
\end{prop}

\begin{proof}(sketch)
We have
\begin{align*}
(F^\star \Omega &F)(j,k) 
= \sum_l (F^\star \Omega)(j,l)F(l,k)\\
&= \frac{1}{n}
\sum_l (\omega^{(j-1)(l-1)}\omega^{l-1}) \bar\omega^{(l-1)(k-1)}
&= 
\frac{1}{n}\sum_{r=0}^{n-1}\omega^{r(j-k+1)} \\
&=
\begin{cases}
1 &\text{if $j\equiv k-1 \mod{n}$} \\
0 &\text{otherwise}.
\end{cases}
\end{align*}
This completes the proof of the first part of the proposition.
 
Since $F$ is unitary, from the first part we get $\pi F^\star = F^\star \Omega$. The second part of the proposition follows from this equation. 
\end{proof}

\begin{prop}
{\bf Spectral decomposition for circulants.} Let
$c:= (c_1,\dots,c_n)$.
Then $\Circ(c) = F^\star p_c(\Omega) F$. 
Hence the eigen-pairs of $\Circ(c)$ are 
$p_c(\omega^{k-1})$, $f_k$,
for $k=1,2,\dots,n$, where $f_k$ is the $k$th column of $F^\star$. 
\end{prop}

\begin{proof}
Recall $\Circ(c) = p_c(\pi)$.
We then have 
$$
\Circ(c) = p_c(\pi) 
= p_c(F^\star \Omega F)
= F^\star p_c(\Omega) F.
$$
Hence
$(\Circ(c))F^\star = F^\star p_c(\Omega)$.
\end{proof}

\defn Let $S$ be an $n\times n$ matrix. Then $S$ is a {\bf skew-circulant matrix} or {\bf skew-circulant} for short, if there exist scalars $a_1,\dots,a_n$, such that $S = \SCirc(a_1,a_2,\dots,a_n)$ where
\begin{equation*}
\SCirc(a_1,a_2,\dots,a_n):=
\begin{pmatrix}
a_1     & a_2   & a_3    & \dots & a_n \\
-a_n     &a_1   &a_2    &\dots &a_{n-1} \\
-a_{n-1} &-a_n   &a_1    &\dots &a_{n-2} \\
\vdots  &\vdots &\vdots &\ddots &\vdots \\
-a_2 &-a_3 &-a_4 &\dots &a_1 
\end{pmatrix}.
\end{equation*}

Note that we can get the skew-circulant $\SCirc(a_1,\dots,a_n)$ from the circulant $\Circ(a_1,\dots,a_n)$ by changing the sign of all elements below the main diagonal. 

Note that the $n\times n$ skew-circulants form a linear subspace of the space of all $n\times n$ matrices.

\defn We use the notation $\eta:=\eta_n:=\SCirc(0,1,0,\dots,0)$ for the {\bf basic skew-circulant} with size $n$. 

Just as the basic circulant $\pi$ plays a fundamental role in the theory of circulants, so the basic skew-circulant $\eta$ plays a fundamental role in the theory of skew-circulants. 
Here is a picture of $\eta$ when $n=4$:
$$
\eta_4 =
\begin{pmatrix}
0     &1   &0    &0 \\
0     &0   &1    &0 \\
0     &0   &0    &1 \\
-1    &0   &0    &0 
\end{pmatrix}.
$$

The following result is easy to verify. I omit the proof.

\begin{prop}
The basic $n\times n$ skew-circulant $\eta$ satisfies the following equations:
$$
\eta^n = -I,\quad \eta^T = \eta^\star = \eta^{-1} = -\eta^{n-1}.
$$
The minimum polynomial of $\eta$ is $\lambda^n+1$.
For every sequence of scalars $a_1,a_2,\dots,a_n$,
$$
\SCirc(a_1,a_2,\dots,a_n)
= a_1 I + a_2 \eta + a_3 \eta^2 + \cdots a_n \eta^{n-1}.
$$
\end{prop}

Note that, if $a:=(a_1,\dots,a_n)$ then
$\SCirc(a) = p_a(\eta)$. From this it is clear that all skew-circulants of the same order commute. In other words, the circulants of a given order form a commutative algebra generated by the single matrix $\eta$. 

\Notation Let $n$ be a positive integer and let $\sigma:=\sigma_n := \exp(i \pi/n) = \cos(\pi/n) + i\sin(\pi/n)$. Let
$\Omega^{1/2}:=\Omega_n^{1/2}:= \Diag(1,\sigma_n,\sigma_n^2,\dots,\sigma_n^{n-1})$.

Note that $\omega_n=\sigma_n^2$ and hence
$\Omega=(\Omega^{1/2})^2$. 

\begin{prop}{\bf The relation between the basic circulant and the basic skew-circulant.}
The basic matrices $\pi$ and $\eta$ satisfy the following equation:
$$ 
\eta=\sigma \Omega^{1/2}\pi \bar{\Omega}^{1/2}
 =\Omega^{1/2}(\sigma\pi) \bar{\Omega}^{1/2}.
$$
\end{prop}

Note that the matrix $\Omega^{1/2}\pi \bar{\Omega}^{1/2}$ is unitarily similar to the matrix $\pi$.
Hence the matrix $\eta$ is a scaled version of a matrix that is unitarily similar to the matrix $\pi$.

\begin{proof}
We have
\begin{align*}
\sigma &\Diag(1,\sigma,\sigma^2,\dots,\sigma^{n-1}) \pi \Diag(1,\bar{\sigma},\bar{\sigma}^2,\dots,\bar{\sigma}^{n-1})\\
&= \sigma
\begin{pmatrix}
0      &\bar{\sigma} &0             &\dots &0 \\
0      & 0            &\bar{\sigma} &\dots &0 \\
0      & 0            &0             &\dots &0 \\
\vdots &\vdots        &\vdots        &      &\vdots\\
0      & 0            &0             &\dots &\bar{\sigma} \\
\sigma^{n-1} & 0      &0             &\dots &0
\end{pmatrix}
=\eta.
\end{align*}
\end{proof}

\Notation Let the $n\times n$ matrix $H:=H_n$, be defined by 
$$H^\star := \Omega^{1/2} F^\star.$$

Here is a picture of $H^\star$:
$$
H^\star := \frac{1}{\sqrt n}
\begin{pmatrix}
1      &1        &1           & \dots 
&1 \\
1      &\sigma\omega   &\sigma\omega^2    &\dots 
&\sigma\omega^{n-1} \\
1      &\sigma^2\omega^2 &\sigma^2\omega^4  &\dots 
&\sigma^2\omega^{2(n-1)} \\
\vdots &\vdots   &\vdots      &\ddots 
&\vdots \\
1      &\sigma^{n-1}\omega^{n-1}&\sigma^{n-1}\omega^{2(n-1)} &\dots 
&\sigma^{n-1}\omega^{(n-1)(n-1)} 
\end{pmatrix}.
$$

\begin{prop}
The matrix $H$ is unitary.
\end{prop}

\begin{proof}
Recall the $F$ is unitary. Note that $\Omega^{1/2}$ is unitary. Hence $H^\star = \Omega^{1/2} F^\star$ is   unitary.
\end{proof}

\begin{prop}{\bf Spectral decomposition of the basic skew-circulant.}
The basic skew-circulant matrix $\eta$ satifies the equation $\eta = H^\star (\sigma\Omega) H$. 
Hence the eigen-pairs of $\eta$ are $\sigma^{2k-1}$, 
$h_k$, for $k=1,\dots,n$, where $h_k$ is the $k$th column of the matrix $H^\star$. 
\end{prop}

\begin{proof}
Using the relation between $\eta$ and $\pi$ and the spectral decomposition of $\pi$, we have
$$
\eta = \sigma\Omega^{1/2}\pi \bar{\Omega}^{1/2}
=\sigma \Omega^{1/2}F^\star \Omega F \bar{\Omega}^{1/2}
=\sigma H^\star \Omega H.
$$
\end{proof}

\begin{prop}
{\bf Spetral decomposition for skew-circulants.} 
Let $a:=(a_1,\dots,a_n)$. Then
$\SCirc(a)=H^\star p_a(\sigma\Omega)H$. 
Hence the eigen-pairs of $\Circ(a)$ are 
$p_a(\sigma^{2k-1})$, $h_k$,
for $k=1,2,\dots,n$, where $h_k$ is the $k$th column of the matrix $H^\star$.
\end{prop}

\begin{proof}
We have
$$
\SCirc(a)=p_a(\eta)=p_a(H^\star(\sigma\Omega)H)
=H^\star p_a(\sigma\Omega)H.
$$
\end{proof}

\section*{References}

\begin{itemize}

\item
Behn, A.; Driessel, K.R.; and Hentzel, I.R.
(2011)
The eigen-problem for some special near-Toeplitz centro-skew tridiagonal matrices,
arXiv:1101.88347v1[mathSP]
\smallskip

\item
Behn, A.; Driessel, K.R.; Hentzel, I.R.; Vander Velden, K.A.; and Wilson, J.
(2011)
Some nilpotent, tridiagonal matrices with a special sign pattern,
Submitted to \emph{Linear Algebra and Appl.}
\smallskip

\item
Catral, M.; Olesky, D.D.; and van den Driesche, P. (2009)
Allow problems concernng spectral properties of sign pattern matrices, 
\emph{Linear Algebra and Appl.}
430, 3080-3094
\smallskip

\item
Davis, P.J.(1979), 
\emph{Circulant Matrices}, 
Wiley
\smallskip

\item
Drew, J.H.; Johnson, C.R.; Olesky, D.D.; and van den Driesche, P. (2000)
Spectrally arbitrary patterns, 
\emph{Linear Algebra and Appl.}
308, 121-137
\smallskip

\item
Elsner, L.; Olesky, D.D.; and van den Driesche, P. (2003)
Low rank perturbations and the spectrum of a tri-diagonal sign pattern,
\emph{Linear Algebra and Appl.}
308, 121-137
\smallskip

\item
Halmos, P.R. (1958) 
\emph{Finite-Dimensional Vector Spaces},
Van Nostrand
\smallskip

\item
Kailath, T.; Kung, S.Y.; and Morf, M. (1979a)
Displacement ranks of matrices and linear equations,
\emph{J. Math. Anal. Appl.}
68, 395-407
\smallskip

\item
Kailath, T.; Kung, S.Y.; and Morf, M. (1979b)
Displacement ranks of a matrix,
\emph{Bull. Amer. Math. Soc.}
1, 769-773
\smallskip

\item
Kailath, T. and Sayed, A.H. (1999)
\emph{Fast Reliable Algorithms for Matrices with Structure},
SIAM
\smallskip

\item
Trench (2004)
Characterization and properties of matrices with generalized symmetry or skew symmetry,
\emph{Lin. Alg. Appl.} 377, 207-218
\smallskip

\end{itemize}

\end{document}